\documentclass[11pt]{article}
\usepackage{mathrsfs}
\usepackage{amsmath}
\usepackage{amssymb}
\usepackage{graphicx}
\usepackage{epic}
\renewcommand{\paragraph}{\roman{paragraph}}

\def \la{\lambda}

\newtheorem{theorem}{\scshape \mdseries  Theorem}[section]
\newtheorem{lemma}[theorem]{\scshape \mdseries  Lemma}
\newtheorem{coro}[theorem]{\scshape \mdseries  Corollary}

\begin{document}

\title{\sf  The least eigenvalue of signless Laplacian of non-bipartite graphs with given domination number\thanks{
Supported by National Natural Science Foundation of China (11071002, 11371028, 71101002), Program for New Century Excellent
Talents in University (NCET-10-0001), Key Project of Chinese Ministry of Education (210091),
Specialized Research Fund for the Doctoral Program of Higher Education (20103401110002),
Project of Educational Department of Anhui Province (KJ2012B040),
Scientific Research Fund for Fostering Distinguished Young Scholars of Anhui University(KJJQ1001).}
}
\author{Yi-Zheng Fan$^{1,}$\thanks{Corresponding author.
 E-mail addresses: fanyz@ahu.edu.cn(Y.-Z. Fan), tansusan1@aiai.edu.cn (Y.-Y. Tan).}, \ Ying-Ying Tan$^2$\\
    {\small  \it $1$. School of Mathematical Sciences, Anhui University, Hefei 230601, P. R. China}\\
  {\small  \it $2$. Department of Mathematics and Physics, Anhui University of Architecture, Hefei 230601, P. R. China} \\
 }
\date{}
\maketitle

\noindent {\bf Abstract:}
Let $G$ be a connected non-bipartite graph on $n$ vertices with domination number $\gamma \le \frac{n+1}{3}$.
We investigate the least eigenvalue of the signless Laplacian of $G$, and present a lower bound for such eigenvalue in terms of
the domination number $\gamma$.



\section{Introduction}

Let $G=(V(G),E(G))$ be a simple graph with vertex set $V(G)=\{v_1,v_2,\ldots,v_n\}$ and edge set $E(G)$.
The {\it adjacency matrix} of $G$ is defined to be
a $(0,1)$-matrix $A(G)=[a_{ij}]$,
where $a_{ij}=1$ if $v_i$ is adjacent to $v_j$, and
$a_{ij}=0$ otherwise.
The {\it degree matrix} of $G$ is defined by $D(G)=\hbox{diag}(d_G(v_1), d_G(v_2),\cdots, d_G(v_n))$, where $d_G(v)$ or simply $d(v)$ is the degree
of a vertex $v$ in $G$.
The matrix $Q(G)=D(G)+A(G)$ is called the {\it signless Laplacian matrix} (or {\it $Q$-matrix}) of $G$.
It is known that $Q(G)$ is nonnegative, symmetric and positive semidefinite, so its eigenvalue are real and can be arranged as:
$q_1(G)\geq q_2(G)\geq \cdots \geq q_n(G)\geq 0.$
We call the eigenvalues of $Q(G)$ as the {\it $Q$-eigenvalues} of $G$ and refer the readers to \cite{cve0,cve1,cve2,cve3,cve4} for the survey on this topic.
The least $Q$-eigenvalue $q_n(G)$ is denoted by $q_{\min}(G)$, and the eigenvectors corresponding to $q_{\min}(G)$ are called the {\it first $Q$-eigenvectors} of $G$.

If $G$ is connected, then $q_{\min}(G)=0$ if and only if $G$ is bipartite.
So, here we are concerned  with the least eigenvalue of connected non-bipartite graph.
 Desai and Rao \cite{des} use the least $Q$-eigenvalue to characterize the bipartiteness of graphs.
As a consequence of this work, Shaun and Fan \cite{shaun} establish the relationship between the least $Q$-eigenvalue and some parameters such as vertex or edge bipartiteness;
and in \cite{fan2} they present some upper bounds for the  least $Q$-eigenvalue in terms of the edge bipartiteness and minimize the signless Laplacian spread.
Cardoso et al. \cite{car} minimize the least $Q$-eigenvalue of non-bipartite graphs.
Liu et al. \cite{liu} minimize the least $Q$-eigenvalue of non-bipartite unicyclic graphs with fixed number of pendant vertices.
Lima et al. \cite{lim} survey the known results and present some new ones for the least $Q$-eigenvalues of graphs.
Our research group investigate how the least $Q$-eigenvalue changes when relocating bipartite branches \cite{wang},
which provides an easier way to get some known results on this topic, and also characterize the unique graph whose least
$Q$-eigenvalue attains the minimum among all non-bipartite graphs with fixed number of pendant vertices \cite{fan}.
In a more general setting, the least eigenvalue of the Laplacian of mixed graph has been discussed in \cite{fan1, tan}.

Recall that a vertex set $S$ of the graph $G$ is called a {\it dominating set} if every vertex of $V(G)\backslash S$ is adjacent to at least one  vertex of $S$.
The {\it domination number} of $G$, denoted by $\gamma(G)$, is the minimum of the cardinalities of all domination sets in $G$.
Surely, if $G$ has no isolated vertices then $\gamma \leq \frac{n}{2}$.
With respect to the adjacency matrix, Stevanovi\'c et al. \cite{stev} determine the unique graph with maximal spectral radius
among all graphs with no isolated vertices and fixed domination number;
Zhu \cite{zhu} characterize the unique graph whose least eigenvalue achieves the minimum among all graphs with fixed domination number.
With respect to the Laplacian matrix,
Brand and Seifter \cite{bra} give a upper bound for the spectral radius in terms of domination number.
Lu et al. \cite{lu}, Nikiforov et al. \cite{niki} and Feng \cite{feng} give some  bounds for the second least eigenvalue 
and the spectral radius of graphs, respectively.
In addition, Feng et al. \cite{feng1} minimize the Laplacian spectral radius among all trees with given domination number.

But few work appears on the relation between the signless Laplacian eigenvalue and the domination number, except that He and Zhou \cite{he} use domination number to upper bound the least signless Laplacian eigenvalue.
In this paper, we will investigate the lower bound for the least $Q$-eigenvalue of a non-bipartite graph in terms of the domination number.
For convenience, a graph is called {\it minimizing} in a certain graph class if its least $Q$-eigenvalue attains the minimum among all graphs in the class.
Denote by $\mathscr{G}^\gamma_n$ the set of all connected non-bipartite graphs of order $n$ with the domination number $\gamma$,
and by $\mathscr{G}^\gamma_n(g)\;(g <n)$ the set of graphs in $\mathscr{G}^\gamma_n$ for which the minimum length of odd cycles is $g$.
When $\gamma \le \frac{n+1}{3}$, we characterize the unique minimizing graph among all graphs in $\mathscr{G}^\gamma_n$,
and hence provide a lower bound for the least $Q$-eigenvalue in terms of the domination number.

\section{Preliminaries}

We first give some preliminary knowledge and notations. Denote by $C_n$ a cycle and $P_n$ a path, both on $n$ vertices.
A graph $G$ is called {\it trivial} if it contains only one vertex;
otherwise, it is called {\it nontrivial}. A graph $G$ is called {\it unicyclic} if it is connected and contains exactly
one cycle. The minimum length of all cycles in $G$ is called the {\it girth} of $G$.
A {\it pendant vertex} of $G$ is a vertex of degree 1 and a {\it quasi-pendant vertex} is one adjacent to a pendant vertex.

Let $x=(x_1, x_2, \ldots, x_n) \in \mathbb{R}^n$.
Denote $|x|=(|x_1|,|x_2|, \ldots, |x_n|)$.
Let $G$ be a graph on vertices $v_1,v_2,\ldots, v_n$. The vector $x$ can be considered
as a function defined on $V(G)$, which maps each vertex $v_i$ of $G$ to the value $x_i$, i.e. $x(v_i)=x_i$.
If $x$ is an eigenvector of $Q(G)$, then it defines on $G$ naturally, i.e. $x(v)$ is the entry of $x$ corresponding to $v$.
One can find
that the quadratic form $x^TQ(G)x$ can be written as
$$x^TQ(G)x=\sum_{uv \in E(G)}[x(u)+x(v)]^2. \eqno(2.1)$$
The eigenvector equation $Q(G)x=\la x$ can be interpreted as
$$[\la -d(v)] x(v)= \sum_{u \in N_G(v)} x(u) \hbox{~ for each~} v \in V(G), \eqno(2.2)$$
where $N_G(v)$ denotes the neighborhood of $v$ in $G$. In addition,
for an arbitrary unit vector $x \in \mathbb{R}^n$,
$$q_{\min} (G)\leq x^TQ(G)x, \eqno(2.3)$$
with equality if and only if $x$ is a first $Q$-eigenvector of $G$.

Let $G_1$, $G_2$ be two vertex-disjoint graphs, and let $v\in V(G_1)$, $u\in V(G_2)$. The {\it coalescence} of
$G_1$ and $G_2$ with respect to $v$ and $u$, denoted by $G_1(v)\diamond G_2(u)$, is obtained from $G_1$ and
$G_2$ by identifying $v$ with $u$ and forming a new vertex $p$, which is also denoted as $G_1(p)\diamond G_2(p)$.
If a connected graph $G$ can be expressed as $G=G_1(p)\diamond G_2(p)$,
where $G_1$ and $G_2$ are nontrivial subgraphs of $G$ both containing $p$,
then $G_1$(or $G_2$) is called a {\it branch} of $G$ with root $p$.
Let $G=G_1(v_2)\diamond G_2(u)$ and $G^*=G_1(v_1)\diamond G_2(u)$, where $v_1$ and $v_2$ are two distinct vertices of $G_1$
and $u$ is a vertex of $G_2$. We say that $G^*$ is obtained from $G$ by {\it relocating $G_2$ from $v_2$ to $v_1$}.


%
%
%
%

\begin{lemma}\label{relocate}
Let $G_1$ be a connected graph containing at least two vertices $v_1$, $v_2$, and let $G_2$ be a connected bipartite graph containing a
vertex $u$. Let $G=G_1(v_2)\diamond G_2(u)$ and $G^*=G_1(v_1)\diamond G_2(u)$. If there exist a first $Q$-eigenvector of $G$ such that
$|x(v_1)|\geq |x(v_2)|$, then
$q_{\min}(G^*) \leq q_{\min}(G),$
with equality only if $|x(v_1)|=|x(v_2)|$, $d_{G_2(u)}x(u)=-\Sigma_{v\in N_{G_2}(u)}x(v)$, and $G^*$ has a first $Q$-eigenvector $\widetilde{x}$ such that
$|\widetilde{x}|=|x|$.
\end{lemma}

{\it Proof:}
By Lemma 2.5 in \cite{wang}, it suffices to prove when the equality holds $G^*$ has a first $Q$-eigenvector $\widetilde{x}$ such that
$|\widetilde{x}|=|x|$.
This can be easily verified if retracing the proof of Lemma 2.5 in \cite{wang}. \hfill $\blacksquare$


\begin{lemma}{\em \cite{wang}} \label{value}
Let $G$ be a connected non-bipartite graph, and let $x$ be a first $Q$-eigenvector of $G$. Let $T$ be a
tree, which is a nonzero branch of $G$ with respect to $x$ and with root $u$. Then $|x(q)| < |x(p)|$ whenever $p, q$
are vertices of $T$ such that $q$ lies on the unique path from $u$ to $p$.
\end{lemma}

\section{Minimizing the least $Q$-eigenvalue among all graphs in $\mathscr{G}^\gamma_n$}
Denote by $U_n^k(g)$ the unicyclic graph of order $n$, which is also obtained from an odd cycle $C_g \;(g<n)$ and a star $S_{1,k}$ by adding a path $P_l$ connecting
(or identifying) one vertex of the cycle and the center of the star (if $l=1$), where $l=n+1-g-k$; see Fig. 3.1.
Surely, if $k \ge 2$, then $$\gamma(U_n^k(g)) \le \gamma(U_n^{k-1}(g)) \le \cdots \le \gamma(U_n^1(g))=:\gamma_g,$$
where $\gamma_g=\lceil\frac{n-1}{3}\rceil$ or $\gamma_g=\lceil\frac{n}{3}\rceil$ depending on whether $3 |g $ or not.

Fixed $n$ and odd $g \in [3,n-1]$, for each $\gamma \in [\lceil\frac{g}{3}\rceil, \gamma_g]$,
there exists one or  more graphs $U_n^k(g)$ with domination number $\gamma$;
the unique one with minimum number $k$ among these graphs is denoted by $V_n^\gamma(g)$.


\begin{center}
\includegraphics[scale=.6]{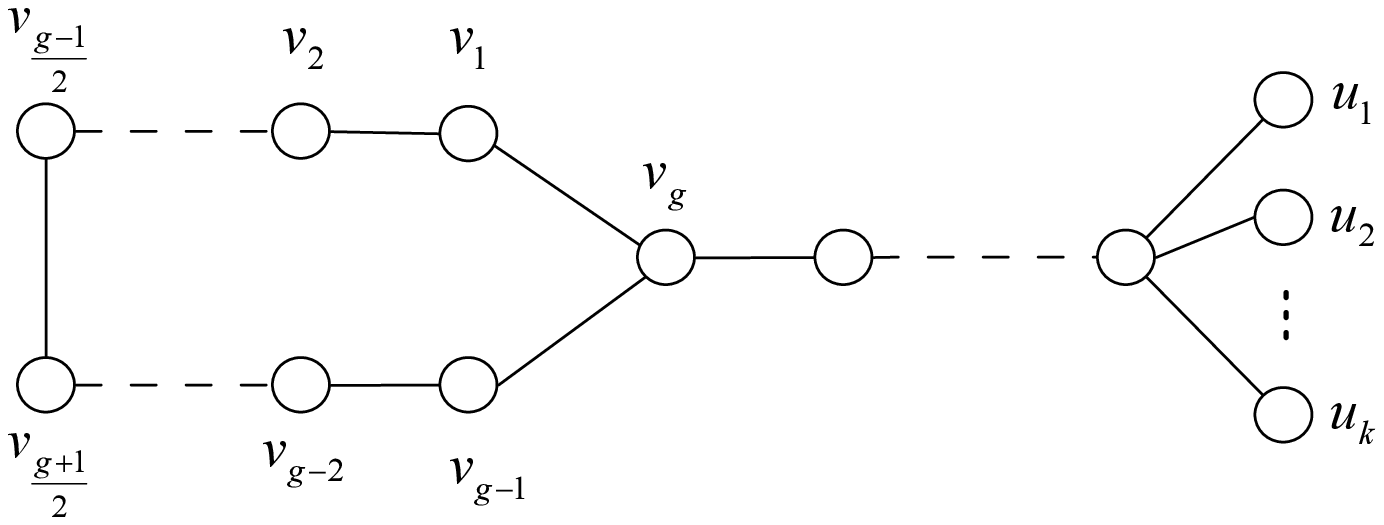}

{\small Fig. 3.1. The graph $U_n^k(g)$ with $g <n$.}
\end{center}

\begin{lemma}{\em\cite{fan}} \label{sign}
Let $U^k_n(g)$ be the graph with some vertices labeled as in Fig. 3.1,
where $v_1, v_2, \cdots, v_g$ are the vertices of the unique cycle $C_g$
labeled in an anticlockwise way. Let $x$ be a first $Q$-eigenvector of $U^k_n(g)$. Then

\noindent{\em (1)} $x(v_i)=x(v_{g-i})$ for $i=1,2, \cdots, \frac{g-1}{2}$;

\noindent{\em (2)} $x(v_{\frac{g-1}{2}})x(v_{\frac{g+1}{2}})> 0$, and $x(v)x(w)< 0$ for other edges $vw$ of $U^k_n(g)$ except $v_{\frac{g-1}{2}}v_{\frac{g+1}{2}}$;

\noindent{\em (3)} $|x(v_g)|> |x(v_1)|> |x(v_2)|> \cdots >|x(v_{\frac{g-1}{2}})|> 0$.

\end{lemma}

Denote by $\mathscr{U}_n^k(g)$ the set of unicyclic graphs of order $n$ with odd girth $g$ and $k \ge 1$ pendant vertices.

\begin{lemma}{\em\cite{fan}} \label{minpengraph}
Among all graphs in $\mathscr{U}_n^k(g)$, $U^k_n(g)$ is the unique minimizing graph.
\end{lemma}

\begin{lemma}{\em\cite{fan}} \label{minpen}
The least $Q$-eigenvalue  of $U^k_n(g)$ is strictly increasing with respect to $k\geq 1$ and odd $g\geq 3$, respectively.
\end{lemma}

\begin{coro} \label{decr-gamma}
The least $Q$-eigenvalue  of $V_n^\gamma(g)$ is strictly decreasing with respect to $\gamma$.
\end{coro}

{\it Proof:}
Suppose that $\gamma \ge \lceil\frac{g}{3}\rceil +1$, $V_n^\gamma(g)=:U_n^k(g)$, $V_n^{\gamma-1}(g)=:U_n^{k'}(g)$.
Clearly, $k < k'$, and the result follows by Lemma \ref{minpen}.
\hfill $\blacksquare$

\begin{coro} \label{uv}
Let $\gamma(U^k_n(g))=\gamma$. If $U^k_n(g) \ne V_n^\gamma(g)$, then $q_{\min}(U^k_n(g))> q_{\min}(V_n^\gamma(g))$.
\end{coro}

{\it Proof:}
Suppose that $V_n^\gamma(g)=U_n^{k'}(g)$.
If $U_n^{k}(g) \ne V_n^{\gamma}(g)$, then $k > k'$, and the result follows by Lemma \ref{minpen}.
\hfill $\blacksquare$

\begin{coro} \label{decr-girth}
The least $Q$-eigenvalue  of $V_n^\gamma(g)$ is strictly increasing with respect to odd $g$.
\end{coro}

{\it Proof:} Suppose that $g \ge 5$.
In the graph  $V_n^\gamma(g)=:U_n^k(g)$ as listed in Fig. 3.1, replacing the edge $v_{g-2}v_{g-1}$ by $v_{g-2}v_1$, we obtain a graph
$G' \in \mathscr{U}_n^{k+1}(g-2)$.
If letting $x$ be a unit first $Q$-eigenvector of $V_n^\gamma(g)$,
then $$q_{\min}(V_n^\gamma(g))=x^TQ(V_n^\gamma(g))x=x^TQ(G')x \ge q_{\min}(G') \ge q_{\min}(U_n^{k+1}(g-2)),$$
where the second equality holds as $x(v_1)=x(v_{g-1})$ by Lemma \ref{sign}, and the last inequality holds by Lemma \ref{minpengraph}.

Furthermore, $q_{\min}(V_n^\gamma(g)) > q_{\min}(G')$; otherwise, $x$ is also a first $Q$-eigenvector of $G'$,
and by considering the eigenvector equations of $V_n^\gamma(g)$ and  $G'$ on the vertex $v_{g-1}$ both associated with $x$,
we will have $x(v_{g-1})=-x(v_{g-2})$; a contradiction to Lemma \ref{sign}.

Notice that $\gamma(U_n^{k+1}(g-2))=:\gamma' \le \gamma(U_n^k(g))=\gamma$.
So, by Corollary \ref{uv} and Corollary \ref{decr-gamma},
$$q_{\min}(V_n^\gamma(g)) > q_{\min}(G') \ge q_{\min}(U_n^{k+1}(g-2)) \ge q_{\min}(V_n^{\gamma'}(g-2)) \ge q_{\min}(V_n^{\gamma}(g-2)).$$
\hfill $\blacksquare$

\begin{lemma} \label{unispan}
Let $G \in \mathscr{G}^\gamma_n$. Then $G$ contains a non-bipartite spanning unicyclic subgraph with domination number $\gamma$.
\end{lemma}

{\it Proof:} If $\gamma=1$, the result is easily verified. So we assume that $\gamma \ge 2$.
Let $U=\{u_1,u_2,\ldots,u_\gamma\}$ be a dominating set of $G$ of size $\gamma$, let $W=V(G) \backslash U$.
Let $B$ be a bipartite spanning subgraph of $G$, which is obtained by deleting all possible edges within $U$ or $W$.

{\bf Case 1:} Suppose that $B$ is connected. Thus, there exist two vertices in $U$, say $u_1$ and $u_2$, such that $N_B(u_1)\cap N_B(u_2)\neq \emptyset$.
Assume that $w_1 \in N_B(u_1)\cap N_B(u_2)$.
Deleting all edges between $u_2$ and the vertices of $(N_B(u_1)\cap N_B(u_2))\backslash \{w_1\}$ (if it is nonempty),
we will get a subgraph $B_1$ of $B$ such that $u_2$ shares exactly one neighbor with $u_1$.
If $U\backslash \{u_1,u_2\} \ne \emptyset$,
noting that $B_1$ is also connected, there exists one vertex $w_2 \in N_B(u_1)\cup N_B(u_2)$ such that $w_2$ is adjacent to one vertex, say $u_3$ in $U\backslash \{u_1,u_2\}$.
Deleting all edges between $u_3$ and  the vertices of $(N_{B_1}(u_3) \backslash  \{w_2\}) \cap (N_B(u_1)\cup N_B(u_2))$,
we will get a subgraph $B_2$ of $B_1$ such that $u_3$ shares exactly one neighbor with exactly one of $u_1$ and $u_2$.
Repeating the above process, we will arrive at a subgraph $B_{\gamma-1}$ of $B$ such that for each $i=2,3, \ldots,\gamma$,
$u_i$ shares exactly one neighbor with exactly one of $u_1, u_2, \ldots, u_{i-1}$.
So $B_{\gamma-1}$ is a tree with domination number $\gamma$.
Since $G$ is non-bipartite, there exists at least one edge $e$ within $U$ or $W$.
Adding the edge $e$ to $B_{\gamma-1}$, the resulting graph is as desired.


{\bf Case 2:} Suppose that $B$ is not connected.
Let $B_1,B_2,\ldots,B_k$ be the components of $B$ with bipartitions  $(U_1, W_1), (U_2, W_2), \ldots, (U_k, W_k)$ respectively.
Since $G$ is connected, there exists a spanning tree $B_1$ of $G$ obtained from $B$ by adding $k-1$ edges between $U_i$s and $U_j$s, or $W_i$s and $W_j$s.
As $G$ is non-bipartite, adding the edges of $E(G)\backslash E(B_1)$ to $B_1$ such that the first odd cycle $C$ appears, we arrive at a graph $B_2$.
If $B_2$ contains only one cycle (i.e. the cycle $C$), the result follows.
Otherwise, $B_2$ contains some even cycles, each of which must contain an edge with endpoints both within $U_i$s or $W_j$s.
Deleting one of such edges will break an even cycle and preserve an odd cycle despite the deleted edge shared by $C$ or not.
Repeating the process of breaking even cycles, we finally arrive at a connected subgraph containing exactly one odd cycle (not necessarily being $C$), as desired.
%
%
%
%
\hfill $\blacksquare$

Denote by $\mathscr{V}_n^\gamma(g)$ the set of unicyclic graphs of order $n$ with odd girth $g<n$ and domination number $\gamma$.

\begin{theorem}\label{minuni}
Among all graph in $\mathscr{V}_n^\gamma(g)$, where $\lceil\frac{g}{3}\rceil \le \gamma \le \gamma_g$, $V_n^\gamma(g)$ is the unique minimizing graph.
\end{theorem}

{\it Proof:}
Let $G$ be a minimizing graph in $\mathscr{V}_n^\gamma(g)$, and let $x$ be a unit first $Q$-eigenvector of $G$.
In order to obtain the result, it suffices to prove $G$ contains exactly one {\it pendant star}
(i.e. the star with maximum possible size centered at a quasi-pendant vertex),
thus $G=U_n^k(g)$ for some $k$ and $U_n^k(g)=V_n^\gamma(g)$ by Corollary \ref{uv}.

Suppose that $G$ contains two or more than two quasi-pendant vertices.
Among all pendant stars of $G$, choose two pendant stars attached at $p$ and $q$ respectively such that one of them has maximum size.
Without loss of generality, assume that $|x(p)| \ge |x(q)|$.
Relocating the pendant star attached at $q$ to $p$, we will get a graph $G_1$ such that $\varrho(G_1) > \varrho(G)$ and $\gamma(G) \ge \gamma(G_1)$,
where $\varrho(G)$ denote the maximum size of the pendant stars of $G$.
By Lemma \ref{relocate}, we also have
  $q_{\min}(G_1)\leq q_{\min}(G)$.

If the graph $G_1$ has more than one quasi-pendant vertices, repeating the above process of relocating pendant stars, we will get a sequence of graphs
$G=G_0,G_1,\ldots,G_m$ such that from $G_{i-1}$ to $G_i$ the pendant star at $q_{i-1}$ is relocated to $p_{i-1}$ for $i=1,2,\ldots,m$, where $p_0:=p$ and $q_0:=q$, and
$$\gamma=\gamma(G)\ge \gamma(G_1) \ge \cdots \ge \gamma(G_m)=:\gamma', \; q_{\min}(G)\ge q_{\min}(G_1) \ge \cdots \ge q_{\min}(G_m),$$
where $G_m$ contains exactly one pendant star, i.e. $G_m=U_n^{k}(g)$ for some $k$.
This can be done as $\varrho(G_i)$ is strictly increasing and is bounded by a finite number.


Now by the above discussion and Corollaries \ref{decr-gamma} and \ref{uv}, we have
$$ q_{\min}(G) \geq q_{\min}(G_{m-1}) \ge q_{\min}(U_n^{k}(g)) \ge q_{\min}(V_n^{\gamma'}(g)) \geq q_{\min}(V_n^\gamma(g)).\eqno(3.1)$$
Since $G$ is the minimizing graph in $\mathscr{V}_n^\gamma(g)$, all the inequalities in (3.1) become equalities.
So by Corollary \ref{decr-gamma} $\gamma=\gamma'$, and by Corollary \ref{uv} $U_n^{k}(g)=V_n^\gamma(g)$.
Also by the equality $ q_{\min}(G_{m-1}) = q_{\min}(U_n^{k}(g))$, if letting $y$ be a first $Q$-eigenvector of $G_{m-1}$,
then by Lemma \ref{relocate}, $|y(p_{m-1})|=|y(q_{m-1})|$, and $U_n^{k}(g)$ has a first $Q$-eigenvector $z$ such that $|z|=|y|$.
So, $|z(p_{m-1})|=|z(q_{m-1})|$.
Note that  $p_{m-1}$ is exactly the quasi-pendant vertex of $U_n^{k}(g)$, and cannot have value (given by $z$) equal to any other vertex by Lemma \ref{sign} and Lemma \ref{value}; a contradiction.
 \hfill $\blacksquare$

\begin{theorem}\label{main-g}
Among all graphs in $\mathscr{G}^\gamma_n(g)$, where $\lceil\frac{g}{3}\rceil \le \gamma \le \gamma_g$,
$V_n^\gamma(g)$ is the unique minimizing graph.
\end{theorem}

{\it Proof:}
Let $G$ be a minimizing graph in $\mathscr{G}^\gamma_n(g)$.
By Lemma \ref{unispan}, $G$ contains a non-bipartite spanning unicyclic subgraph with domination number $\gamma$; denoted by $\widetilde{G}$.
If $\widetilde{G}=C_n$, noting that $G \ne C_n$ by the definition of $\mathscr{G}^\gamma_n(g)$, then there exists an edge $uv \in E(G)\backslash E(\widetilde{G})$ joining two vertices of $C_n$.
Then $C_n$ is split into two cycles $C^1,C^2$ by the edge $uv$, where $C^1$ is odd and $C^2$ is even.
Noting that $\gamma(G)=\gamma(C_n+uv)=\gamma(C_n)=\gamma$, there exists one edge $u'v' \in E(C^2)\backslash E(C^1)$ such that
$\gamma(C_n+uv-u'v')=\gamma$.

So we assume $\widetilde{G}$ contains pendant vertices with girth $\widetilde{g}$, where $g \le \widetilde{g} <n$.
By Theorem \ref{minuni} and Corollary \ref{decr-girth},
$$ q_{\min}(G) \ge q_{\min}(\widetilde{G}) \ge q_{\min}(V_n^\gamma(\widetilde{g})) \ge q_{\min}(V_n^\gamma(g)).$$
If $ q_{\min}(G)=q_{\min}(V_n^\gamma(g))$, then $\widetilde{g}=g$ by  Corollary \ref{decr-girth}, and $\widetilde{G}=V_n^\gamma(\widetilde{g})=V_n^\gamma(g)$ by Theorem \ref{minuni}.
 %
%
Returning to the origin graph $G$, which is obtained from $\widetilde{G}=V_n^\gamma(g)$ possibly by adding some edges.
Assume that $E(G)\backslash E(V_n^\gamma(g))\neq \emptyset$.
Let $x$ be a unit first $Q$-eigenvector of $G$. Then
\begin{align*}
q_{\min}(G) &= \sum_{uv \in E(G)}[x(u)+x(v)]^2\\
& =\sum_{uv \in E(V_n^\gamma(g))}[x(u)+x(v)]^2+\sum_{uv \in E(G)\backslash E(V_n^\gamma(g))}[x(u)+x(v)]^2\\
& \geq \sum_{uv \in E(V_n^\gamma(g))}[x(u)+x(v)]^2 \geq q_{\min}(V_n^\gamma(g)).
\end{align*}
Since $q_{\min}(G)=q_{\min}(V_n^\gamma(g))$, $x$ is also the first $Q$-eigenvector of $V_n^\gamma(g)$, and $x(u)+x(v)=0$ for each edge
$uv \in E(G)\backslash E(V_n^\gamma(g))$, which cannot occur by Lemma \ref{sign} and Lemma \ref{value}.
So, $G=V_n^\gamma(g)$ is the unique minimizing graph.
\hfill $\blacksquare$

\begin{coro}
Among all graphs in $\mathscr{G}^\gamma_n$, where $1\leq \gamma\leq \frac{n+1}{3}$,
$V_n^\gamma(3)$ is the unique minimizing graph.
\end{coro}

{\it Proof:}
Let $G$ be a minimizing graph in $\mathscr{G}^\gamma_n$.
Assume that the minimum length of odd cycles of $G$ is $g$.
If $g=n$, then $G=C_n$; but the graph $U_n^1(n-2)$ has the same domination number as $C_n$, and
$q_{\min}(C_n) > q_{\min}(U_n^1(n-2))$ by  a similar discussion as in the proof of Corollary \ref{decr-girth}; a contradiction.
So, we may assume $g <n$.
By Theorem \ref{main-g} and Corollary \ref{decr-girth}
$$ q_{\min}(G) \ge q_{\min}(V_n^\gamma(g)) \ge q_{\min}(V_n^\gamma(3)).$$
Since $G$ is minimizing, all the inequalities become equalities, so $G=V_n^\gamma(g)$ by Theorem \ref{main-g} and
$g=3$ by Corollary \ref{decr-girth}, which implies that $G=V_n^\gamma(3)$.
\hfill $\blacksquare$

If $n \ge 3\gamma+1$, then $V_n^\gamma(3)=U_n^{n-3\gamma}(3)$.
If $n=3\gamma-1$ or $n=3\gamma$, then $V_n^\gamma(3)=U_n^1(3)$, that is, $V_n^\gamma(3)$ is obtained from $C_3$ by appending a path $P_{n-2}$.

We finally arrive at the main result of this paper.

\begin{coro}
Let $G$ be a connected non-bipartite graph of order $n$ with domination number $\gamma \leq \frac{n+1}{3}$.
Then $q_{\min}(G) \ge q_{\min}(V_n^\gamma(3))$, with equality if and only if $G=V_n^\gamma(3)$.
\end{coro}


\small

\begin{thebibliography}{90}
\bibitem{bra} C. Brand, N. Seifter, Eigenvalues and domination in graphs, {\it Math. Slovaca}, 46(1)(1996), 33-39.

\bibitem{car} D. M. Cardoso, D. Cvetkovi\'{c}, P. Rowlinson, S. K. Simi\'{c}, A sharp lower bound for the least eigenvalue of the signless
Laplacian of a non-bipartite graph, {\it Linear Algebra Appl.}, 429(2008), 2770-2780.

\bibitem{cve0} D. Cvetkovi\'{c}, P. Rowlinson, S. K. Simi\'{c}, Eigenvalue bounds for the signless Laplacian,
{\it Publ. Inst. Math.} (Beograd), 81(95)(2007), 11-27.

\bibitem{cve1} D. Cvetkovi\'{c}, P. Rowlinson, S. K. Simi\'{c}, Signless Laplacians of finite graphs,
{\it Linear Algebra Appl.}, 423(2007), 155-171.

\bibitem{cve2} D. Cvetkovi\'{c}, S. K. Simi\'{c}, Towards a spectral theory of graphs based on the signless Laplacian, I,
{\it Publ. Inst. Math.} (Beograd), 85(99)(2009), 19-33.

\bibitem{cve3} D. Cvetkovi\'{c}, S. K. Simi\'{c}, Towards a spectral theory of graphs based on the signless Laplacian, II,
{\it Linear Algebra Appl.}, 432(2010), 2257-2272.

\bibitem{cve4} D. Cvetkovi\'{c}, S. K. Simi\'{c}, Towards a spectral theory of graphs based on the signless Laplacian, III,
{\it Appl. Anal. Discrete Math.}, 4(2010), 156-166.

\bibitem{des} M. Desai, V. Rao, A characterization of the smallest eigenvalue of a graph, {\it J. Graph Theory}, 18(1994), 181-194.

\bibitem{shaun} S. Fallat, Y.-Z. Fan, Bipartiteness and the least eigenvalue of signless Laplacian of graphs, {\it Linear Algebra Appl.}, 436(2012), 3254-3267.

\bibitem{fan} Y.-Z. Fan, Y. Wang, H. Guo, The least eigenvalues of the signless Laplacian of non-bipartite graphs with pendant vertices, {\it Discrete Math.}, 313(2013), 903-909.


\bibitem{fan1} Y.-Z. Fan, S.-C. Gong, Y. Wang, Y.-B. Gao, First eigenvalue and first eigenvectors of a nonsingular unicyclic mixed graph, {\it Discrete Math.}, 309(2009), 2479-2487.


\bibitem{fan2} Y.-Z. Fan, Shaun Fallat, Edge bipartiteness and signless Laplacian spread of graphs, {\it Appl. Anal. Discrete Math.}, 6(2012), 31-45.


\bibitem{feng} L.-H. Feng, G.-H. Yu, X.-Q. Li Laplacian eigenvalues of graphs with given domination number, {\it Asian-European J. Math.}, 2(1)(2009), 71-76.

\bibitem{feng1} L.-H. Feng, G.-H. Yu, Q. Li, Minimizing the Laplacian eigenvalues for trees with given domination number,
{\it Linear Algebra Appl.}, 419(2006), 648-655.


\bibitem{he} C.-X. He, M. Zhou,  A sharp upper bound on the least signless Laplacian eigenvalue using domination number, {\it Graph Combin.}, doi: 10.1007/s00373-013-1330-z.



\bibitem{lim} L.S. de Lima, C.S. Oliveira, N.M.M.de Abreu, V. Nikiforov, The smallest eigenvalue of the signless Laplacian,
{\it Linear Algebra Appl.}, 435(2011), 2570-2584.



\bibitem{liu} R.-F. Liu, H.-X. Wan, J.-J. Yuan, H.-C. Jia,
The least eigenvalue of the signless Laplacian of non-bipartiteness unicyclic graphs with $k$ pendant vertices,
{\it Electron. J. Linear Algebra},  26(2013), 333-344.

\bibitem{lu} M. Lu, H.-Q. Liu, F. Tian, Bounds of Laplacian spectrum of graphs based
on the domination number, {\it Linear Algebra Appl.}, 402(2005) 390-396.

\bibitem{niki} V. Nikiforov, Bounds on graph eigenvalues I, {\it Linear Algebra Appl.}, 420(2007), 667-671.

\bibitem{stev} D. Stevanovi\'c, M. Aouchiche, P. Hansen,
On the spectral radius of graphs with a given domination number, {\it Linear Algebra Appl.}, 428(2008), 1854-1864.

\bibitem{tan} Y.-Y. Tan, Y.-Z. Fan, On edge singularity and eigenvectors of mixed graphs, {\it Acta Math. Sinica Eng. Ser. B}, 24(1)(2008), 139-146.

\bibitem{wang} Y. Wang, Y.-Z. Fan, The least eigenvalues of signless Laplacian of graphs
under perturbation, {\it Linear Algebra Appl.}, 436(2012), 2084-2092.


\bibitem{zhu} B.-X. Zhu, The least eigenvalue of a graph with a given domination number, {\it Linear Algebra Appl.}, 437(2012), 2713-2718.


\end{thebibliography}
\end{document}